\documentclass{amsart}
\usepackage{amsfonts}

\newtheorem{theorem}{Theorem}
\theoremstyle{plain}

\newtheorem{lemma}{Lemma}

\numberwithin{equation}{section}

\begin{document}
\title[On a new exponential bound for M-acceptable variables]{A note on a new exponential bound for M-acceptable random variables}
\author{Cheikhna Hamallah Ndiaye}
\email{ham111266@yahoo.fr}
\address{LERSTAD, Universit\'{e} de Saint-Louis.}
\address{LMA - Laboratoire de Math\'ematiques Appliqu\'ees \\
Universit\'e Cheikh Anta Diop BP 5005 Dakar-Fann S\'en\'egal}
\author{Gane Samb LO}
\email{ganesamblo@ganesamblo.net, ganesamblo@yahoo.com}
\address{Corresponding author.\\
LSTA, UPMC, FRANCE and LERSTAD, Universit\'{e} de Saint-Louis, SENEGAL.\\
Tel +221 33 961 23 40. Fax +221 961 53 53. BP 234 SAINT-LOUIS SENEGAL}

\begin{abstract} \textbf{English : } 
We present a new exponential inequality as a generalization of that of Sung
\textit{et al.} \cite{sun2011} for $M$-acceptable random variables, and hence for
extended negative ones. Our result is based on the simple real
inequality $e^{x} \leq 1+x+(|x|/2)e^{|x|}, x\in\mathbb{R}$, in place of the
following one: $e^{x} \leq 1+x+(x^{2}/2)e^{|x|}, x\in\mathbb{R}$, used by many authors. We compare the given bound with former ones.\\

\noindent \textbf{French : } :
Nous g\'en\'eralisation la notion d'acceptabilit\'e d'une variable al\'eatoire \`a celle de $M$-acceptabilit\'e, puis \'etablissons
une g\'en\'eralisation de l'in\'egalit\'e exponentielle de Sung \textit{et al.} \cite{sun2011} pour ce type de variables, contenant aussi les variables
\'etend\^ument n\'egatives. Notre r\'esultat se base sur l'in\'egalit\'e simple suivante $e^{x} \leq 1+x+(|x|/2)e^{|x|}, x\in\mathbb{R}$, \`a la place de celle-ci :  $e^{x} \leq 1+x+(x^{2}/2)e^{|x|}, x\in\mathbb{R}$, utilis\'ee par plusieurs auteurs. Nous comparons notre borne exponentielle avec d'autres
 pr\'ec\'edemment introduites.\\

\end{abstract}

\keywords{ Exponential inequality, Convergence rate, Almost sure
convergence, $M$-Acceptable random variables, Negatively associated random
variables, Negatively dependent random variables, Extended Negatively
dependent random variables, Laplace transform.}
\subjclass[2010]{60F15, 62G20}
\maketitle

\Large

\section{Introduction}

\label{sec1}

Let \ $\{X_{n},n\geq 1\}$ \ be a sequence of random variables defined on a
fixed probability space $(\Omega ,\mathcal{F},\mathbb{P}).$ An exponential
inequality for the partial sum \ ${\sum_{i=1}^{n}}(X_{i}-\mathbb{E}X_{i})$ \
was introduced and proved by Sung \textit{et al.} \cite{sun2011} for a class of dependent
random variables, that is for $M$-\textit{acceptable} random variables. This
type of inequality provides a measure of the convergence rate for the strong
law of large numbers. It is rather classical for the independent case but
seems very recent in the dependent case.\\

\bigskip

\noindent We refer to Sung \textit{et al.} \cite{sun2009} for an extended introduction on
dependent random variables. We will only recall basic definitions to make
this paper self-contained. First, recall that a finite family of random
variables \ $\{X_{i},1\leq i\leq n\}$ is positiveley associated if for any
real functions $f_{1}$ and $f_{2}$\ defined on $\mathbb{R}^{n}$ and coordinatewise
increasing, we have%
\begin{equation*}
Cov(f_{1}(X_{1},...,X_{n}),\ f_{2}(X_{1},...,X_{n}))\geq 0.
\end{equation*}%
This type of dependence was introduced in Esary, Proschan and Walkup \cite{esa}.
Next the negative association was set up by Alam and Saxena \cite{ala} and
carefully studied by Joag-Dev and Proschan \cite{joa}. We have negative
association whenever for any pair of disjoint subsets \ $A$ \ and $B$ \ of \
$\{1,2,...,n\},$ with respective cardinalities $p$ and $q$, and for any pair
functions respectively defined of $\mathbb{R}^{p\text{ }}$and $\mathbb{R}%
^{q}$, we have
\begin{equation*}
Cov(f_{1}(X_{i},i\in A),\ f_{2}(X_{j},j\in B))\leq 0.
\end{equation*}%
Further, the concept of extended negatively dependent random variables was
proposed by Liu \cite{liu} as a generalization of negative dependence introduced
by Lehmann \cite{leh}, meaning that the family\ $\{X_{1},...,X_{n}\}$ \ is
extended negatively dependent if there is some \ $M>0$ \ such that for any
real numbers $x_{1},...,x_{n}$, we have :
\begin{equation*}
\mathbb{P(}X_{1}\leq x_{1},...,X_{n}\leq x_{n})\leq M\ {\prod_{i=1}^{n}}%
\mathbb{P(}X_{i}\leq x_{i})
\end{equation*}%
and
\begin{equation*}
\mathbb{P(}X_{1}>x_{1},...,X_{n}>x_{n})\leq M{\prod_{i=1}^{n}}\mathbb{P(}%
X_{i}>x_{i}).
\end{equation*}%
The special case $M=1$ corresponds to the Lehman original definition.
Finally, Sung \textit{et al.} \cite{sun2011} defined another negative dependence
denamed \textit{acceptable} for \ a family $\{X_{1},...,X_{n}\}$ \ for which
there exists \ $\delta >0$ \ such that for any real \ $\lambda $ \ such that
\ $|\lambda |\leq \delta ,$ \
\begin{equation}
\mathbb{E}\ exp(\lambda {\sum_{i=1}^{n}}X_{i})\leq {\prod_{i=1}^{n}}\mathbb{%
E}\ exp(\lambda X_{i}).  \label{accepta}
\end{equation}

\noindent Before we go further, we recall that from the steps of Joag - Dev and Proschan
\cite{joa}, that a number of well known multivariate distributions fulfils the
negative association property, such as the coordinates of a multinomial
vector, that of a multivariate hypergeometric. Also the following stochastic
processes satisfy the negative dependence :  Dirichlet; permutation
distribution; negatively correlated normal distributions; random sampling
without replacement and joint distribution of ranks; etc.\\

\bigskip

\noindent Now as to the probability inequality field, as pointed out by Sung \textit{et
al.} \cite{sun2011}, exponential inequalities for positively associated random
variables were obtained by Devroye \cite{dev}, Ioannides and Roussas \cite{ioa},
Oliveira \cite{oli}, Sung \cite{sun2007}), Xing and Yang (2008) ; Xing, Yang and Liu \cite{xinyl2008}. We may also cite the following other contributors in this field: Kim and Kim \cite{kim}, Nooghabi and Azarnoosh \cite{noo}, Roussas \cite{rou}, Sung \cite{sun2009}, Xing \cite{xin2009}, Xing and Yang \cite{xin2010}, Xing \textit{et al.} \cite{xinyl2009}. Many of these inequalities are generalized by Sung\textit{et al.} \cite{sun2011}, particularly for negatively associated random variables.\newline

\noindent In this note we aim at improving significally the current exponential
inequalities into two directions. The first is to extend the acceptability of
Sung \textit{et al.}\cite{sun2011} by $M$-\textit{acceptability, for }$M>0,$ by replacing (%
\ref{accepta}) by%
\begin{equation}
\mathbb{E}\ exp(\lambda {\sum_{i=1}^{n}}X_{i})\leq M{\prod_{i=1}^{n}}\mathbb{%
E}\ exp(\lambda X_{i}).  \label{acceptam}
\end{equation}%
It is clear that extending these inequalities to $M$-\textit{acceptable }%
random variable may result in a significant generalization as we point out
below. The second generalization concerns the inequality itself. We will see
in the next section that our results are based on a much precise real analysis
inequality. For example Sung \textit{et al.} \cite{sun2011} used the inequality
\begin{equation}
\forall (x\in R),e^{x}\leq 1+x+\frac{x^{2}}{2}e^{\left\vert x\right\vert }.
\label{ineq1}
\end{equation}%
To obtain a better exponential inequality for negative random variables, we
will use this inequality%
\begin{equation}
\forall (x\in R),e^{x}\leq 1+x+\frac{\left\vert x\right\vert }{2}%
e^{\left\vert x\right\vert }.  \label{ineq2}
\end{equation}

\noindent The remainder of this paper is organized as follows. In Section %
\ref{sec2}, we state our results and make some comments on their generality.
The proofs, including that of (\ref{ineq2}), are given in Section \ref{sec3}. It is worth mentioning that most of the methods of 
Sung \textit{et al.} \cite{sun2011} are reconducted here and hence, will not detailed as it should be.

\bigskip

\section{Results and comments}

\label{sec2}

\subsection{Results}

\label{sec21}

As pointed out in the introduction, the inequality (\ref{ineq2}) is
essential. It is used to prove the

\begin{lemma} \label{lemma}
Let \ $X$ \ be a random variable with \ $\mathbb{E}e^{\delta |X|}<\infty $ \
for some \ $\delta >0.$ Then for any \ $0<\lambda \leq \delta /2,$
\begin{equation*}
\mathbb{E}e^{\lambda (X-\mathbb{E}X)}\leq exp(K\lambda )
\end{equation*}%
where \ $K$ \ is defined as \ $K=(\mathbb{E}|X|^{2})^{1/2}\ \mathbb{E}\
e^{\delta |X|}.$
\end{lemma}

\noindent This lemma is itself a generalization of Lemma 1 in Sung \textit{et al.} \cite{sun2011}
(2011). \bigskip We will get the new exponential inequality from it later.

\begin{theorem} \label{theo1}
Let $\{X_n, \ n \geq 1\}$ be a sequence of identically distributed $M$
acceptable random variable with $\mathbb{E} \ e^{\lambda|X_1|} < \infty$ for
some $\delta > 0$.

Then for any $\varepsilon \geq K$,
\begin{equation*}
\mathbb{P(}\left\vert \sum_{i=1}^{n}\ (X_{i}-\mathbb{E}\ X_{i})\right\vert
>n\ \varepsilon )\leq 2\ M\ exp(-n\ \delta /2(\varepsilon -K)),
\end{equation*}%
where $K=(\mathbb{E}|X_{1}|^{2})^{1/2}\ \mathbb{E}\ e^{\delta |X_{1}|}$.
\end{theorem}

\subsection{Comments}

\label{sec22}

\subsubsection{Comparison with Sung \textit{et al.} results}

For the special case of negatively associated random variables, let us
compare our result with those of Sung et al. \cite{sun2011}. In their Remark 2.2, Sung \textit{et al.} \cite{sun2011} pointed out that Sung \cite{sun2009} upper bound and their own one, are of the same convergence rate when $\varepsilon$ is bounded. Our results results based on an unbounded $\varepsilon$ is hence more general than those in the two papers. Next, for $\varepsilon_n = 2 \ (K \alpha \ log n/n)^{1/2}$, $\alpha > 0$ and $K = 2(E|X_1|^2)^{1/2} \ E \ e^{\delta|X_1|}$, they obtained the bound $2 exp(-\alpha \ log n)$ in their Theorem 2.2. They remarked its supiority to those of Kim and Kim \cite{kim}, Nooghabi and Azarnook \cite{noo}, Xing \cite{xin2009}, Xing et al. \cite{xinyl2009} and Xing and Yang \cite{xin2010}. It is clear that our upper bound is better than that of 
\textit{Sung et al.} \cite{sun2011} and then to all the cited ones above.

\subsubsection{Particular case}

By choosing $\varepsilon = n^{\alpha -1} + K$ in Theorem \ref{theo1}, we
have the following result

\begin{theorem} \label{theo2}
\label{theo2} Let $\{X_{n},\ n\geq 1\}$ be a sequence of identically
distributed $M$ acceptable random variables with $\mathbb{E}\ e^{\delta
|X_{1}|}<\infty $ for some $\delta >0$. Let $\varepsilon _{n}=n^{\alpha -1}+K
$, where $\alpha >0$ and $K=(\mathbb{E}|X_{1}|^{2})^{1/2}\ \mathbb{E}\
e^{\delta |X_{1}|}$. Then for all $n$,
\begin{equation*}
\mathbb{P}\bigl(\bigr|\sum_{i=1}^{n}\ (X_{i}-\mathbb{E}\ X_{i})\bigr|>n\
\varepsilon _{n})\bigr)
\end{equation*}

\begin{equation*}
\leq 2 \ M \ exp\bigl(-\frac{\delta}{2} \ n^\alpha\bigr)
\end{equation*}
\end{theorem}

\noindent \textbf{Proof}. Let $\varepsilon _{n}=n^{\alpha -1}+K$, where $%
\alpha >0$ and $K=(E|X_{1}|^{2})^{1/2}\ E\ e^{\delta |X_{1}|}$. Then $%
\varepsilon _{n}-K\geq 0$ for all $n$. Hence the result follows from Theorem \ref{theo1}.

\section{Proofs}

\label{sec3}

Let us begin with the proof of Inequality (\ref{ineq2}).

\bigskip

\subsection{Proof of (\ref{ineq2})} We first remark that for all $x\in R,$%

\begin{equation}
\frac{|x|}{2}\ e^{|x|}=\frac{|x|}{2}+\frac{|x|^{2}}{2}+{\sum_{n=2}^{\infty }}%
\frac{|x|^{n+1}}{2n!}.  \label{ineq3}
\end{equation}%
Next, we remark that for $n\geq 3,$ $2(n-1)!\leq n(n-1)!=n!$, and then $%
1/n!\leq 1/(2(n-1)!)$ and next%
\begin{equation*}
\frac{x^{n}}{n!}\leq \frac{|x|^{n}}{n!}\leq \frac{|x|^{n}}{2(n-1)!}
\end{equation*}%
and finally%
\begin{equation}
{\sum_{n=3}^{\infty }}\frac{x^{n}}{n!}\leq {\sum_{n=3}^{\infty }}\frac{%
|x|^{n}}{2(n-1)!}={\sum_{n=3}^{\infty }}\frac{|x|^{n}}{2(n-1)!}.
\label{ineq4}
\end{equation}%
To conclude, we use the decomposition%
\begin{equation*}
e^{x}=1+x+\frac{\left\vert x\right\vert ^{2}}{2}+{\sum_{n=3}^{\infty }}\frac{%
|x|^{n}}{2(n-1)!}
\end{equation*}%
and use (\ref{ineq3}) and (\ref{ineq4}) to get%
\begin{equation*}
e^{x}=1+x+\frac{\left\vert x\right\vert ^{2}}{2}+{\sum_{n=3}^{\infty }}\frac{%
x^{n}}{n!}\leq e^{x}\leq 1+x+\frac{\left\vert x\right\vert ^{2}}{2}+{%
\sum_{n=3}^{\infty }}\frac{|x|^{n}}{2(n-1)!}
\end{equation*}%
\begin{equation*}
=1+x+\frac{\left\vert x\right\vert }{2}(\left\vert x\right\vert +{%
\sum_{n=3}^{\infty }}\frac{|x|^{n-1}}{(n-1)!})=1+x+\frac{\left\vert
x\right\vert }{2}(\left\vert x\right\vert +{\sum_{n=2}^{\infty }}\frac{%
|x|^{n}}{n!})
\end{equation*}%
\begin{equation*}
=1+x+\frac{\left\vert x\right\vert }{2}({\sum_{n=1}^{\infty }}\frac{|x|^{n}}{%
n!})\leq 1+x+\frac{\left\vert x\right\vert }{2}({\sum_{n=0}^{\infty }}\frac{%
|x|^{n}}{n!})=1+x+\frac{\left\vert x\right\vert }{2}e^{\left\vert
x\right\vert }.
\end{equation*}

\noindent Now let us turn to the proof of Lemma \ref{lemma}.

\bigskip

\subsection{Proof of Lemma \ref{lemma}}
In view of inequality (\ref{ineq2}) and by using the H\"{o}lder inequality for $0<\lambda \leq
\delta /2,$ to get%
\begin{equation*}
\mathbb{E}e^{\lambda (X-\mathbb{E}X)}\leq 1+\lambda \mathbb{E}(X-\mathbb{E}%
X)+\frac{\lambda }{2}\mathbb{E}(|X-\mathbb{E}X|)\ e^{\lambda |X-\mathbb{E}%
X|})
\end{equation*}%
\begin{equation*}
=1+\frac{\lambda }{2}\mathbb{E}(|X-\mathbb{E}X|\ e^{\lambda |X-\mathbb{E}X|})
\end{equation*}%
\begin{equation*}
\leq 1+\frac{\lambda }{2}(\mathbb{E}(X-\mathbb{E}X)^{2})^{1/2}(\mathbb{E}%
e^{2\lambda |X-\mathbb{E}X|})^{1/2}.
\end{equation*}%
Now%
\begin{equation*}
(\mathbb{E}e^{2\lambda |X-\mathbb{E}X|})^{1/2}\leq (\mathbb{E}e^{2\lambda
|X+\left\vert \mathbb{E}X\right\vert })^{1/2}\leq (\mathbb{E}e^{2\lambda
\left\vert X\right\vert })^{1/2}(e^{2\lambda \left\vert \mathbb{E}%
X\right\vert })^{1/2}
\end{equation*}%
and by Jensen inequality, $e^{2\lambda \left\vert \mathbb{E}X\right\vert
}\leq \mathbb{E}e^{2\lambda \left\vert X\right\vert }$, we end up with
\begin{equation*}
(\mathbb{E}e^{2\lambda |X-\mathbb{E}X|})^{1/2}\leq \mathbb{E}e^{2\lambda
\left\vert X\right\vert }.
\end{equation*}%
Next,%
\begin{equation*}
\mathbb{E}(X-\mathbb{E}X)^{2})^{1/2}\leq (\mathbb{E}X^{2}-(\mathbb{E}%
X)^{2})^{1/2}\leq (\mathbb{E}X^{2})^{1/2}.
\end{equation*}%
Putting all this together, we get%
\begin{equation*}
\mathbb{E}e^{\lambda (X-\mathbb{E}X)}\leq 1+\frac{\lambda }{2}(\mathbb{E}%
X^{2})^{1/2}\mathbb{E}e^{2\lambda \left\vert X\right\vert }\leq 1+\frac{%
\lambda }{2}K\leq \exp (\frac{\lambda }{2}K).
\end{equation*}%

\noindent This achieves the proof. Let us close the proofs section by that of Theorem \ref{theo1}.

\subsection{Proof of Theorem \ref{theo1}} Let $\varepsilon \geq K$. By
Markov's inequality, by the definition of $M$ acceptable random variables
and by the Lemma, we have for $0<\lambda \leq \delta /2$,%
\begin{equation*}
\mathbb{P}(\sum_{i=1}^{n}\ (X_{i}-\mathbb{E}\ X_{i})>n\ \varepsilon )=%
\mathbb{P}(exp(\lambda \sum_{i=1}^{n}\ (X_{i}-\mathbb{E}\
X_{i}))>exp(\lambda \ n\ \varepsilon ))
\end{equation*}%
\begin{equation*}
\leq exp(-\lambda \ n\ \varepsilon )\ E\ exp(\lambda \ \sum_{i=1}^{n}\
(X_{i}-\mathbb{E}\ X_{i}))
\end{equation*}%
\begin{equation*}
\leq M\ exp(-\lambda \ n\ \varepsilon )\ \prod_{i=1}^{n}\ \mathbb{E}%
(exp(\lambda (X_{i}-\mathbb{E}\ X_{i}))\leq M\ exp(-\lambda \ n\ \varepsilon
)\ \prod_{i=1}^{n}\ exp(K\lambda ),
\end{equation*}%
since the $X_{i}^{{}}$ have the same first and second moments. We finally get%
\begin{equation*}
\mathbb{P}(\sum_{i=1}^{n}\ (X_{i}-\mathbb{E}\ X_{i})>n\ \varepsilon )\leq M\
exp(-\lambda \ n\ \varepsilon +\lambda K\ n).
\end{equation*}%
The right member is minimum for $\lambda =\delta /2$ since $\varepsilon \geq
K.$ We then get for this value%
\begin{equation*}
\mathbb{P}(\sum_{i=1}^{n}\ (X_{i}-\mathbb{E}\ X_{1})>n\ \varepsilon )\leq M\
exp(-n\ \delta /2\ (\varepsilon -K)).
\end{equation*}%
Since $\{-X_{n}\ ,\ n\geq 1\}$ are also $M$ acceptable random variables, we
can replace $X_{i}$ by $-X_{i}$ in the above statement that is ,
\begin{equation}
\mathbb{P(}-\sum_{i=1}^{n}\ (X_{i}-\mathbb{E}\ X_{i})n\ \varepsilon )\leq M\
exp(-n\ \delta /2\ (\varepsilon -K)).  \label{l}
\end{equation}

\noindent Now observing that $\left\vert x\right\vert =\max (x,-x)$ and $(\left\vert
x\right\vert >a)=(x>a$ $or$ $x<a)$, we arrive at
\begin{equation*}
\mathbb{P(}\left\vert \sum_{i=1}^{n}\ (X_{i}-\mathbb{E}\ X_{1})\right\vert
>n\ \varepsilon )
\end{equation*}%
\begin{equation*}
\leq \mathbb{P(}\sum_{i=1}^{n}\ (X_{i}-\mathbb{E}\ X_{n})>n\ \varepsilon )+%
\mathbb{P(-}\sum_{i=1}^{n}\ (X_{i}-\mathbb{E}\ X_{n})>n\ \varepsilon )
\end{equation*}%
\begin{equation*}
\leq 2M\ exp(-n\ \delta /2\ (\varepsilon -K)).
\end{equation*}


\begin{thebibliography}{99}
\bibitem{ala} Alam, K. and Saxena, K.M.L.(1981) Positive dependence in
multivariate distributions. \textit{Communications in Statistics Theory and
Methods}, 10, 1183 - 1196.


\bibitem{dev} Dovroye, L.(1991) Exponential inequalities in
non-parametric estimation, In G. Roussas (Ed.) : \textit{Nonparametric functional
estimation and related topics}, (pp. 31 - 44). Dordrecht : Kluwer Academic
Publishers.

\bibitem{esa} Esary, J. D., Proschan, F. and Walkup, D. W.(1967)
Association of random variables, with applications. \textit{Annals of Mathematical
Statistics}, \textbf{(38)}, 1466 - 1474.

\bibitem{ioa} Ioannides, D. A. and Roussas, G. G.(1999) Exponential
inequality for associated random variables. \textit{Statistics and Probability
Letters}, \textbf{(42)}, 423 - 431.

\bibitem{joa} Joag - Der, K. and Proschan, F.(1983). Negative association of
random variables with applications. \textit{Ann. of Statist.}, \textbf{(11)}, 286 - 295.

\bibitem{kim} Kim, T. S. and Kim, H. C.(2007) On the exponential inequality
for negative dependent sequence. \textit{Communications of Korean Mathematical
Society}, \textbf{(22)}, 315 - 321.

\bibitem{leh} Lehmann, E.(1966) Some concepts of dependence. \textit{Annals of
Mathematical Statistics}, \textbf{(37)}, 1137 - 1153.

\bibitem{liu} Liu, L.(2009) Precise large deviatime for dependent random
variables with heany tails. \textit{Statist. Probab. Lett.}, \textbf{(79)}, $n^{o}$ 9, 1290 - 1298.

\bibitem{noo} Nooghabi, H. J. and Azarnoosh, H. A.(2009). Exponential
inequality for negatively associated random variables. \textit{Statistical Papers},
\textbf{(50)}, 419 - 428.

\bibitem{oli} Oliveira, P. E.(2005) An exponential inequality for
associated variables. \textit{Statistics and Probability Letters}, \textbf{(73)}, 189 - 197.

\bibitem{rou} Roussas, G. G.(1996). \textit{Exponential probability inequalities
with some applications}. IMS Lecture Notes Monograph Series, 30, 303 - 319.

\bibitem{sun2007} Sung, S. H.(2007) A note on the exponential inequality for
associated random variables. \textit{Statistics and Probability Letters}, \textbf{(77)}, 1730 -
1736.

\bibitem{sun2009} Sung, S. H.(2009) An exponential inequality for negatively
random variables. \textit{Journal of Inequalities and Applications}, Article ID 64927,
7 pages.

\bibitem{sun2011} Sung, S. H., Srisuradetchai, P. and Volodin, A.(2011) A note
on the exponential inequality for a class of dependent random variables.
\textit{Journal of the Korean Statistical Society}, \textbf{(40)}, 109 - 114.

\bibitem{xin2009} Xing, G.(2009) On the exponential inequalities for strictly
stationary and negatively associated random variables. \textit{Journal of
statistical Planning and inference}, \textbf{(139)}, 3453 - 3460.

\bibitem{xin2008} Xing, G. and Yang, S.(2008) Notes on the exponential
inequalities for strictly and positively associated random variables.
\textit{Journal of Statistical Planning and inference}, \textbf{(138)}, 4132 - 4140)

\bibitem{xin2010} Xing, G. and Yang, S.(2010). An exponential inequalities for
strictly stationnary and negatively associated random variables.
\textit{Communications in Statistics Theory and Methods}, \textbf{(39)}, 340 - 349.

\bibitem{xinyl2008} Xing, G., Yang, S. and Lin, A.(2008). Exponential
inequalities for positively associated random variables and applications.
\textit{Journal of Inequalities and Applications}. Article ID 385362 11 pages.

\bibitem{xinyl2009} XIng, G. Yang, S. Lin, A. and Wang, X.(2009) A remark on the
exponential inequality for negatively associated random variables. \textit{Journal
of the Korean Statistical Society}, \textbf{(38)}, 53 - 57.

\end{thebibliography}
\end{document}